\documentclass[12pt,reqno]{amsart}

\usepackage[arrow,matrix,curve]{xy}

\usepackage[dvips]{graphicx} 

\usepackage{amssymb, latexsym, amsmath, amscd, array,
}

\theoremstyle{definition}

\numberwithin{equation}{section}

\newcommand\N {{\mathbb N}} 

\newcommand\R {{\mathbb R}}
\newcommand\Q {{\mathbb Q}}

\newcommand\RRR{\text{I\!I\!R}}

\begin{document}

\thispagestyle{empty}

\title[Stevin numbers and reality]{Stevin numbers and reality}

\author{Karin Usadi Katz and Mikhail G. Katz$^{0}$}

\address{Department of Mathematics, Bar Ilan University, Ramat Gan
52900 Israel} \email{katzmik@macs.biu.ac.il}




\bigskip\bigskip\bigskip\noindent
\begin{abstract}
We explore the potential of Simon Stevin's numbers, obscured by
shifting foundational biases and by 19th century developments in the
arithmetisation of analysis.
\end{abstract}

\maketitle

\tableofcontents

\section{From discrete arithmetic to arithmetic of the continuum}

Simon Stevin (1548-1620) initiated a systematic approach to decimal
representation of measuring numbers, marking a transition from a
discrete arithmetic as practiced by the Greeks, to the arithmetic of
the continuum taken for granted today, see A.~Malet \cite{Mal06} and
Naets \cite{Nae}.

For over two centuries now, such numbers have been called {\em
real\/}.  Concerns about the {\em reality\/} of numbers generally
preoccupy cognitive scientists and philosophers more than
mathematicians.  Thus, cognitive scientists view mathematical infinity
as necessarily a metaphor~\cite{LN}, while philosophers such as
G.~Hellman~\cite{He} have attempted nominalistic reconstructions that
seek to diminish an investigator's reliance on ontological assumptions
that provoke tensions with philosophical examinations of the
foundations.  Such reconstructions have, in turn, been criticized by
other philosophers \cite{BR}, see also \cite{Ch} for a response.  It
is interesting to note in this context that C. S. Peirce thought of
the Weierstrassian doctrine of the limit as a nominalistic
reconstruction, see J.~Dauben~\cite{Da82} and Section~\ref{DP} below.%
\footnote{\label{great}C.~Boyer refers to Cantor, Dedekind, and
Weierstrass as ``the great triumvirate'', see \cite[p.~298]{Boy}.  The
triumvirate reconstruction of analysis as a nominalistic project is
explored in our text \cite{KK11a}.}

Many mathematicians regard such ontological questions as of limited
relevance to the practice of mathematics.  They feel that a
mathematician reasons the same way, whether or not he thinks
mathematical objects actually pass any reality check, if such were
possible.  At the same time, they readily admit serious negative
effects in the past caused by an undue influence of a preoccupation
with whether such things as complex numbers really {\em exist\/}.

Our goal here is neither to pursue the cognitive thread, nor to
endorse any nominalistic reconstruction, but rather to focus on the
reception of Stevin's ideas, and how such reception was influenced by
received notions of what a continuum should, or rather should not, be.
We also examine the effects of Platonist perceptions of the real
numbers on the practice of both mathematics and the history of
mathematics, as well as the attitude toward infinitesimal-enriched
extensions of the traditional number system.  Some related issues are
analyzed by Katz and Tall in \cite{KT}.

\section{Stevin's construction of the real numbers}
\label{stevin}

Stevin created the basis for modern decimal notation in his 1585 work
{\em De Thiende\/} (``the art of tenths'').  He argued that quantities
such as square roots, irrational numbers, surds, negative numbers,
etc., should all be treated as numbers and not distinguished as being
different in nature. He wrote that ``there are no absurd, irrational,
irregular, inexplicable or surd numbers.''  He further commented as
follows:
\begin{quote}
It is a very common thing amongst authors of arithmetics to treat
numbers like~$\sqrt{8}$ and similar ones, which they call absurd,
irrational, irregular, inexplicable or surds etc and which we deny to
be the case for number which turns up.
\end{quote}
Thus, Stevin explicitly states that numbers that are {\em not\/}
rational have equal rights of citizenship with those that are.
According to van der Waerden, Stevin's
\begin{quote}
general notion of a real number was accepted, tacitly or explicitly,
by all later scientists \cite[p.~69]{van}.
\end{quote}
D.~Fearnley-Sander wrote that 
\begin{quote}
the modern concept of real number [...] was essentially achieved by
Simon Stevin, around 1600, and was thoroughly assimilated into
mathematics in the following two centuries \cite[p.~809]{Fea}.
\end{quote}
D.~Fowler points out that
\begin{quote}
Stevin [...] was a thorough-going arithmetizer: he published, in 1585,
the first popularization of decimal fractions in the West [...]; in
1594, he described an algorithm for finding the decimal expansion of
the root of any polynomial, the same algorithm we find later in
Cauchy's proof of the intermediate value theorem \cite[p.~733]{Fo}.
\end{quote}

The algorithm is discussed in more detail in \cite[\S10,
p.~475-476]{Ste}.  Unlike Cauchy, who {\em halves\/} the interval at
each step (see Section~\ref{cau}), Stevin subdivides the interval into
{\em ten\/} equal parts, resulting in a gain of a new decimal digit of
the solution at every iteration of the algorithm.%
\footnote{Stevin's numbers were anticipated by E. Bonfils in 1350, see
S.~Gandz \cite{Ga}.  Bonfils says that ``the unit is divided into ten
parts which are called Primes, and each Prime is divided into ten
parts which are called Seconds, and so on into infinity''
\cite[p.~39]{Ga}.}
Thus, while Cauchy's algorithm can be described as a binary search,
Stevin's approach is a more general divide-and-conquer algorithm.

Fowler makes the following additional points (see \cite{Fo}).  The
belief that all arithmetic operations, as well as extracting roots,
etc., should follow the ``same'' rules as the rationals, originates
precisely with Stevin.  The rigorous justification of such a belief
had to await Dedekind's contribution at the end of the 19th century.
The ``existence'' of multiplication of the real numbers was first
proved by Dedekind.  The widespread belief that there exists an {\em
algorithm\/} for determining the digits of the result of multiplying
real numbers in terms of finite pieces of the decimal string, is
unfounded (namely, there {\em is\/} no such algorithm).%
\footnote{To illustrate the point, consider a computer multiplying a
decimal~$.333...$ by $3$, where we are deliberately vague about what
the ellipsis stands for.  A computer programmer can spend an
arbitrarily large time thinking that the resulting decimal will start
with a long string of 9s.  It suffices for a single digit greater than
$3$ to appear at the trillionth rank to prove it wrong, showing that
in this calculation, at no time can the programmer be sure of any
given digit.}
This thread concerning the precise nature of Dedekind's contribution
is pursued further in Section~\ref{DP}.

In his {\em L'Arithm\'etique\/}, Stevin grounds a transition from the
classical arithmetic of the discrete, to a continuous arithmetic, by
means of his well-known ``water-and-wetness" metaphor.  Numbers are
measures, and measures of continuous magnitudes are by their nature
continuous:

\begin{quote}
as well as to continuous water corresponds a continuous wetness, so to
a continuous magnitude corresponds a continuous number" (Stevin, 1585,
see \cite[p.~3]{St85}; quoted in Malet \cite{Mal06}).
\end{quote}

\section{A Stevin-Cauchy proof of the intermediate value theorem}
\label{cau}

How are we to understand van der Waerden's contention that Stevin
numbers were accepted by all later scientists?  To illustrate the
issue at stake, consider Cauchy's proof of the intermediate value
theorem~\cite{Ca21} (the proof was forshadowed in a text of Stevin's,
see Section~\ref{stevin}).  Cauchy constructs an increasing
sequence~$a_n$ and a decreasing sequence~$b_n$ of successive
approximations,~$a_n$ and~$b_n$ becoming successively closer than any
positive distance.%
\footnote{Cauchy's notation for the two sequences is $x_0, x_1, x_2,
\ldots$ and $X, X', X'', \ldots$ \cite[p.~462]{Ca21}.}
At this stage, the desired point is considered to have been exhibited,
by Cauchy.  A modern mathematician may object that Cauchy {\em has
not, and could not have, proved the existence of the limit point\/}.

But imagine that the perspicacious {\em polytechnicien\/} Auguste
Comte%
\footnote{Comte's notes of Cauchy's lectures are extant, see
\cite[p.~437]{Sch}.}
had asked {\em M.~le Professeur\/} Cauchy the following question:
\begin{quote}
Consider a decimal rank, say~$k>0$.  What is happening to the~$k$-th
decimal digit~$a_n^k$ of $a_n$, and~$b_n^k$ of~$b_n$?  
\end{quote}
{\em M.~le Professeur\/} would have either sent Comte to the library
to read Simon Stevin, or else provided a brief argument to show that
for~$n$ sufficiently large, the~$k$-th digit stabilizes, noting that
special care needs to be taken in the case when~$a_n$ is developing a
tail of~$9$s and~$b_n$ is developing a tail of~$0$s.  Clearly the
arguments appearing in Cauchy's textbook are sufficient to identify
the Stevin decimal expression of the limit point.%
\footnote{Note that Cauchy exploited decimal notation on occasion;
see, for instance, \cite[p.~34]{Ca53}.}

From the modern viewpoint, the only item missing is the remark that a
Stevin decimal {\em is\/} a number, {\em by definition\/} (modulo the
technical detail of the identification of the pair of tails).

In the same spirit, D.~Laugwitz points out that in France after 1870,
\begin{quote}
the main objective of French mathematics professors, under the
leadership of the Ecole Polytechnique ... was to prepare students of
engineering and the sciences for useful jobs ... With decimal
expansions of real numbers at hand, nobody was bothered by theories of
irrational numbers \cite[p.~274]{Lau00}.
\end{quote}
The incoherence of triumvirate scholarship in relation to Cauchy was
already analyzed in 1973 by Hourya Benis Sinaceur \cite{Si}, and is
further analyzed in \cite{BK, KK11b} (see also Section~\ref{rival}).

Recall that Stevin's algorithm involved partitioning the interval into
ten parts, and produces an additional rank of the decimal expansion of
the solution at each step of the iteration (see Section~\ref{stevin}).
Much has been said about the proof of the ``existence'' of the real
numbers by the great triumvirate at the end of the 19th century.  But
who needs such an {\em existence\/} proof when the Stevin-Cauchy
method gives an {\em algorithm\/} that produces a concrete infinite
decimal string?

\section{Peirce's framework}

The customary set-theoretic framework (e.g., the language of
equivalence classes of Cauchy sequences; Dedekind cuts; etc.) has
become the reflexive litmus test of {\em mathematical rigor\/} in most
fields of modern mathematics (with the possible exception of the field
of mathematical logic).  Such a framework makes it difficult to
analyze Cauchy's contribution to the foundations of analysis,
particularly Cauchy's use of the concept of an infinitesimal, and to
evaluate its significance.  We will therefore use a conceptual
framework proposed by C.~S.~Peirce in 1897 (going back to his text
{\em How to make our ideas clear\/} of 1878, see
\cite[item~(5.402)]{Pe78}), in the context of his analysis of the
concept of continuity and continuum, which, as he felt at the time, is
composed of infinitesimal parts, see \cite[p.~103]{Ha}.  Peirce
identified three stages in creating a novel concept:

\begin{quote}
there are three grades of clearness in our apprehensions of the
meanings of words.  The first consists in the connexion of the word
with familiar experience. . . . The second grade consists in the
abstract definition, depending upon an analysis of just what it is
that makes the word applicable. . . . The third grade of clearness
consists in such a representation of the idea that fruitful reasoning
can be made to turn upon it, and that it can be applied to the
resolution of difficult practical problems \cite{Pei} (see
\cite[p.~87]{Ha}).
\end{quote}
The ``three grades'' can therefore be summarized as
\begin{enumerate}
\item
familiarity through experience;
\item
abstract definition with an eye to future applications;
\item
fruitful reasoning ``made to turn'' upon it, with applications.
\end{enumerate}

A related taxonomy was developed by D.~Tall \cite{Ta08}, in terms of
his {\em three worlds of mathematics\/}, based on embodiment,
symbolism and formalism in which mathematical proof develops in each
world in terms of recognition, description, definition and deduction.

To apply Peirce's framework to Cauchy's concept of an infinitesimal,
we note that the perceptual stage (1) is captured in Cauchy's
description of continuity of a function in terms of ``varying by
imperceptible degrees''.  Such a turn of phrase occurs both in his
letter to Coriolis of 1837, and in his 1853 text \cite[p.~35]{Ca53}.%
\footnote{Note that both Cauchy's original French ``par degr\'es
insensibles'', and its correct English translation ``by imperceptible
degrees'', are etymologically related to {\em sensory perception\/}.}

At stage (2), Cauchy describes infinitesimals as generated by null
sequences (see \cite{Br}), and defines continuity in terms of an
infinitesimal~$x$-increment resulting in an infinitesimal change
in~$y$.

Finally, at stage (3), Cauchy fruitfully applies the crystallized
concept of an infinitesimal both in Fourier analysis and in evaluation
of singular integrals.  Thus, Cauchy exploits a ``Dirac'' delta
function defined in terms of what would be called today the Cauchy
distribution with an infinitesimal scaling parameter, see
Cauchy~\cite[p.~188]{Ca27}, \cite{Ca27a},
Freudenthal~\cite[p.~136]{Fr}, and Laugwitz~\cite[p.~219]{Lau89}
and~\cite{Lau92}.

Peirce's flexible framework allows us to appreciate Cauchy's
foundational contributions and their fruitful application.  How do
Cauchy's infinitesimals fare in a triumvirate framework?  This issue
is explored in the Section~\ref{straw}.

\section{A case study in triumvirate strawmanship}
\label{straw}

Cauchy's 1821 {\em Cours d'Analyse\/} \cite{Ca21} presented only a
theory of infinitesimals of polynomial rate of growth as compared to a
given ``base'' infinitesimal~$\alpha$.  The shortcoming of such a
theory is its limited flexibility.  Since Cauchy only considers
infinitesimals behaving as polynomials of a fixed infinitesimal,
called the ``base" infinitesimal in 1823, his framework imposes
obvious limitations on what can be done with such infinitesimals.
Thus, one typically can't extract the square root of such a
``polynomial'' infinitesimal.

What is remarkable is that Cauchy did develop a theory to overcome
this shortcoming.  Cauchy's astounding theory of infinitesimals of
arbitrary order (not necessarily integer) is analyzed by Laugwitz
\cite[p.~271]{Lau87}.

In 1823, and particularly in 1829, Cauchy develops a more flexible
theory, where an infinitesimal is represented by an arbitrary {\em
function\/} (rather than merely a polynomial) of a base infinitesimal,
denoted ``$i$".  This is done in Cauchy's 1829 textbook \cite[Chapter
6]{Ca29}.  The title of the chapter is significant.  Indeed, the title
refers to the {\em functions\/} as ``representing" the infinitesimals;
more precisely, {\em ``fonctions qui repr\'esentent des quantit\'es
infiniment petites"\/}.  Here is what Cauchy has to say in~1829:
\begin{quote}
Designons par~$a$ un nombre constant, rationnel ou irrationnel; par
$i$ une quantite infiniment petite, et par~$r$ un nombre variable.
Dans le systeme de quantit\'es infiniment petites dont~$i$ sera la
base, une fonction de~$i$ represent\'ee par~$f(i)$ sera un infiniment
petit de l'ordre~$a$, si la limite du rapport~$f(i)/i^r$ est nulle
pour toutes les valeurs de~$r$ plus petite que~$a$, et infinie pour
toutes les valeurs de~$r$ plus grandes que~$a$ \cite[p. 281]{Ca29}.
\end{quote}

Laugwitz \cite[p.~271]{Lau87} explains this to mean that the order~$a$
of the infinitesimal~$f(i)$ is the uniquely determined real number
(possibly~$+\infty$, as with the function~$e^{-1/t^2}$) such
that~$f(i)/i^r$ is infinitesimal for~$r < a$ and infinitely large
for~$r > a$.

Laugwitz \cite[p.~272]{Lau87} notes that Cauchy provides an example of
functions defined on positive reals that represent infinitesimals of
orders~$\infty$ and~$0$, namely 
\[
e^{-1/i} \mbox{ \ and \ } \frac{1}{\log i}
\]
(see Cauchy \cite[p.~326-327]{Ca29}).

Note that according to P.~Ehrlich's detailed 2006 study \cite{Eh06},
the development of non-Archimedean systems based on orders of growth
was pursued in earnest at the end of the 19th century by such authors
as Stolz and du Bois-Reymond.  These systems appear to have an
anticedent in Cauchy's theory of infinitesimals as developed in his
texts dating from 1823 and~1829.  Indeed, already in 1966, A.~Robinson
pointed out that
\begin{quote}
Following Cauchy's idea that an infinitely small or infinitely large
quantity is associated with the behavior of a function $f(x)$, as $x$
tends to a finite value or to infinity, du Bois-Raymond produced an
elaborate theory of orders of magnitude for the asymptotic behavior of
functions ... Stolz tried to develop also a theory of arithmetical
operations for such entities \cite[p.~277-278]{Ro66}.
\end{quote}
Robinson traces the chain of influences further, in the following
terms:
\begin{quote}
It seems likely that Skolem's idea to represent infinitely large
natural numbers by number-theoretic functions which tend to infinity
(Skolem [1934]), also is related to the earlier ideas of Cauchy and du
Bois-Raymond \cite[p.~278]{Ro66}.
\end{quote}
The reference is to Skolem's 1934 work \cite{Sk}.  

The material presented in the present section, including a detailed
discussion of the chain of influences from Cauchy via Stolz and du
Bois-Reymond and Skolem to Robinson, is contained in an article
entitled ``Who gave you the Cauchy-Weierstrass tale?  The dual history
of rigorous calculus''.  At no point did the article claim that
Cauchy's approach is a variant of Robinson's approach.  Indeed, such a
claim would be preposterous, as Cauchy was not in the possession of
the mathematical tools required to either formulate or justify the
ultrapower construction, requiring as it does a set-theoretic
framework (dating from the end of the 19th century) together with the
existence of ultrafilters (not proved until 1930 by Tarski
\cite{Tar}).

The article was submitted to the periodical ``Revue d'histoire des
sciences" on 5 october 2010.  The article was rejected by editor
Michel Blay five months later, in a letter dated 11 march 2011.  Blay
based his decision on two referee reports.  Referee 1 summarized the
article as follows in his {\em third\/} sentence:
\begin{quote}
Our author interprets A. Cauchy's approach as a formation of the idea
of an infinitely small - a variant of the approach which was developed
in the XXth century in the framework of the nonstandard analysis (a
hyperreal version of E. Hewitt, J. Los, A. Robinson).
\end{quote}
Based on such a strawman version of the article's conception, the
referee came to the following conclusion:
\begin{quote}
From my point of view the author's arguments to support this
conception are quite unconvincing.
\end{quote}
The author indeed finds such such a strawman conception unconvincing,
but the conception was the referee's, not the author's.  The referee
concluded as follows:
\begin{quote}
The fact that the actual infinitesimals lived somewhere in the
consciousness of A. Cauchy (as in many another mathematicians of XIXth
- XXth centuries as, for example , N.N. Luzin) does not abolish his
(and theirs) constant aspiration to dislodge them in the
subconsciousness and to found the calculus on the theory of limit.
\end{quote}
The notion of a Cauchy as a pre-Weierstrassian, apparently espoused by
the referee, is just as preposterous as the notion of Cauchy as a
pre-Robinsonian.  Such a notion is a reflection of a commitment to a
triumvirate ideology, elevated to the status of a conditioned reflex.
Felix Klein knew better: fifty years before Robinson, he clearly
realized the potency of the infinitesimal approach to the foundations
(see Section~\ref{rival}).

Cauchy did not aspire to dislodge infinitesimals; on the contrary, he
used them with increasing frequency in his work, including his 1853
article \cite{Ca53} where he relies on infinitesimals to express the
property of uniform convergence.

The pdf version of the submitted article ``Who gave you, etc.'', as
well as the two referee reports, may be found at the following web
page: http://u.cs.biu.ac.il/$\sim$katzmik/straw.html

Similarly, in his 2007 anthology \cite{Haw}, S.~Hawking reproduces
Cauchy's {\em infinitesimal\/} definition of continuity on page
639--but claims {\em on the same page\/}, in a comic {\em
non-sequitur\/}, that Cauchy ``was particularly concerned to banish
infinitesimals''.

\section{Weierstrassian epsilontics}

If we are to take at face value van der Waerden's evaluation of the
significance of Stevin numbers, what is, then, the nature of
Weierstrass's contribution?  After describing the formalisation of the
real continuum usually associated with the names of Cantor, Dedekind,
and Weierstrass on pages 127-128 of his retiring presidential address
in 1902, E.~Hobson remarks triumphantly as follows:
\begin{quote}
It should be observed that the criterion for the convergence of an
aggregate%
\footnote{i.e. an equivalence class defining a real number}
is of such a character that no use is made in it of infinitesimals
\cite[p.~128]{Ho}.
\end{quote}
Hobson reiterates:
\begin{quote}
In all such proofs [of convergence] the only statements made are as to
relations of finite numbers, no such entities as infinitesimals being
recognized or employed.  Such is the essence of the~$\epsilon$[$,
\delta$] proofs with which we are familiar \cite[p.~128]{Ho}.
\end{quote}
The tenor of Hobson's remarks, is that Weierstrass's primary
accomplishment was the elimination of infinitesimals from foundational
discourse in analysis.  If our students are being dressed to perform
multiple-quantifier epsilontic logical stunts on the pretense of being
taught infinitesimal calculus, it is because infinitesimals are
assumed to be either metaphysically dubious or logically unsound, see
D.~Sherry \cite{She87}.

The significance of the developments in the foundations of the real
numbers at the end of the 19th century was a rigorous proof of the
existence of arithmetic operations, as discussed at the end of
Section~\ref{stevin}.

\section{Dedekind and Peirce}
\label{DP}

The first use of the term {\em real\/} to describe Stevin numbers
seems to date back to Descartes, who distinguished between real and
imaginary roots of polynomials.  Thus, the term was used as a way of
contrasting what were thought of, ever since Stevin, as measuring
numbers, on the one hand, and imaginary ones, on the other.  Gradually
the meaning of the term {\em real number\/} has shifted, to a point
where today it is used in the sense of ``genuine, objective, true
number''.  How legitimate is such usage?

Let us consider what natural science tells us about the physical line.
Descending below the threshold of sensory perception, quantum
physicists tell us that lengths smaller than~$10^{-30}$ meters are not
accessible, even theoretically, to any, existing or future, physical
electron microscope.  This is because the minutest entities considered
even theoretically, such as strings in Witten's M-theory, are never
smaller than a barrier of~$10^{-30}$ (up to a few orders of
magnitude).  Thus, the infinite divisibility taken for granted in the
case of the real line, only holds within a suitable range, even in
principle, in the spatial line.%
\footnote{To put it another way, physical theories such as quantum
mechanics testify to a graininess of physical matter, that is at odds
with the infinite divisibility postulated as a key property of the
real axis.  See also Moore's discussion of Dedekind's position,
below.}
Such tiny real numbers, similarly to infinitesimal quantities, appear
too small to see.  M. Moore%
\footnote{See \cite[p.~82]{Moo07} as well as \cite{Moo02}.}
notes that the minute size of strings
\begin{quote}
does not disqualify strings from admission to our scientific ontology;
and . . . the difference, on this score, between things as small as
strings and distances infinitely small, is one of degree and not of
kind.
\end{quote}

Historians Eves and Newsom make the following claim in Dedekind's
name:
\begin{quote}
Dedekind perceived that {\em the essence of continuity\/} of a
straight line lies in the property that if all the points of the line
are divided into two classes, such that every point in the first class
lies to the left of every point in the second class, then there exists
one and only one point of the line which produces this severance of
the line into the two classes \cite[p.~222]{EN} [emphasis
added--authors].
\end{quote}
However, their claim is at variance with the fact that Dedekind
himself specifically downplayed his claims as to such {\em essence\/},
in the following terms: 
\begin{quote}
If space has at all a real existence, it is {\em not\/} necessary for
it to be [complete] . . . if we knew for certain that space was
[incomplete], there would be nothing to prevent us, in case we so
desired, from filling up its gaps, in thought, and thus making it
[complete] (as cited in \cite[p.~73]{Moo07}).
\end{quote}
Dedekind maintains that his completeness principle gives the essence
of continuity, but denies that we can {\em know\/} that space is
continuous in that sense.  M.~Moore points out that Dedekind
\begin{quote}
shows commendably little sympathy for the idea that we know by
intuition that the line is complete \cite[p.~73]{Moo07}.
\end{quote}
Moore further points out that Dedekind admits that his is ``utterly
unable to adduce a proof of [his account's] correctness, nor has
anyone the power'' (as quoted in \cite[p.~73]{Moo07}).

At variance with the great triumvirate of Cantor, Dedekind, and
Weierstrass,%
\footnote{See footnote~\ref{great}.}
American philosopher Charles Sanders Peirce felt that a construction
of a true continuum necessarily involves infinitesimals.  He wrote as
follows:
\begin{quote}
But I now define a {\em pseudo-continuum\/} as that which modern
writers on the theory of functions call a continuum.  But this is
fully represented by [...]  the totality of real values, rational and
irrational%
\footnote{See CP 6.176, 1903 marginal note.  Here (and below) CP x.y
stands for Collected Papers of Charles Sanders Peirce, volume x,
paragraph y.}
[emphasis added---authors]
\end{quote}
He publicly used the word ``pseudo-continua" to describe real numbers
in the syllabus (CP 1.185) of his lectures on Topics of Logic.  Thus,
Peirce's intuition of the continuum corresponded to a type of a
\mbox{B-continuum} (see Section~\ref{rival}), whereas an A-continuum
to him was a pseudo-continuum.

While Peirce thought of a continuum as being made up of infinitesimal
increments, other authors pursuing B-continuum foundational models
(see Section~\ref{rival}) thought of a real number $x$ as having a
cluster of infinitesimals around it, more precisely a cluster of
points infinitely close to $x$, i.e. differing from $x$ by an
infinitesimal amount.  The alternative to Dedekind's view that a cut
on the rationals corresponds to a {\em single\/} number, is to view
such a cut as being defined by a cluster of infinitely close numbers.

Peirce had a theory of infinitesimals that in many ways anticipated
20th century developments, see J.~Dauben \cite{Da}.  Havenel \cite{Ha}
argues that Peirce's conception was closer to Lawvere's approach than
to Robinson's.  This is corroborated by Peirce's opposition to both
the law of excluded middle, and to a view of the continuum as being
reducible to points.  Both of these points are borne out by the
category-theoretic framework of Lawvere's theory (see J. Bell
\cite{Bel08}), in the context of intuitionistic logic.  On the other
hand, Peirce does not seem to have anticipated the notion of a
nilsquare infinitesimal, which had been anticipated already by
Nieuwentijdt, and implemented in Lawvere's theory, see J.~Bell
\cite{Bel} for details.

\section{Oh numbers, numbers so real}

In an era where no implementation of an infinitesimal-enriched
continuum as yet existed, the successful implementation of real
analysis on the basis of the real number system and~$\epsilon, \delta$
arguments went hand-in-hand with an attempt to ban the infinitesimals,
thought of as an intellectual embarrassment at least since George
Berkeley's time~\cite{Be}.  Hobson is explicit in measuring the
significance of Weierstrass's contribution by the yardstick of the
elimination of infinitesimals.%
\footnote{Today, the didactic value of infinitesimals is becoming
increasingly evident, see \cite{El, KK1, KK2, BH, DM, MD}.  On the
foundational side, an implementation of an infinitisimal-enriched
number system can be presented in a traditional set-theoretic
framework by means, for example, of the ultrapower construction, see
e.g., Keisler~\cite{Ke}; Goldblatt \cite{Go}; M.~Davis~\cite{Dav77}.
Some philosophical implications are explored by B\l aszczyk
\cite{Bl}.}

Today we can perhaps appreciate more clearly, not Weierstrass's, but
Stevin's contribution toward the implementation of quantities that may
have been called Stevin numbers, and that generally go under the
reassuring name of numbers so {\em real\/}.

\section{Rival continua}
\label{rival}

The historical roots of infinitesimals go back to Cauchy, Leibniz, and
ultimately to Archimedes.  Cauchy's approach to infinitesimals is not
a variant of the hyperreals.  Rather, Cauchy's work on the rates of
growth of functions anticipates the work of late 19th century
investigators such as Stolz, du Bois-Reymond, Veronese, Levi-Civita,
Dehn, and others, who developed non-Archimedean number systems against
virulent opposition from Cantor, Russell, and others, see Ehrlich
\cite{Eh06} and Katz and Katz \cite{KK11a} for details.  The work on
non-Archimedean systems motivated the work of T.~Skolem on
non-standard models of arithmetic \cite{Sk}, which stimulated later
work culminating in the hyperreals of Hewitt, \L os, and Robinson.

\begin{figure}
\[
\xymatrix@C=95pt{{} \ar@{-}[rr] \ar@{-}@<-0.5pt>[rr]
\ar@{-}@<0.5pt>[rr] & {} \ar@{->}[d]^{\hbox{st}} & \hbox{\quad
B-continuum} \\ {} \ar@{-}[rr] & {} & \hbox{\quad A-continuum} }
\]
\caption{\textsf{Thick-to-thin: taking standard part (the thickness of
the top line is merely conventional)}}
\label{31}
\end{figure}

Having outlined the developments in real analysis associated with
Weierstrass and his followers, Felix Klein pointed out in 1908 that
\begin{quote}
The scientific mathematics of today is built upon the series of
developments which we have been outlining.  But {\em an essentially
different conception of infinitesimal calculus has been running
parallel with this [conception] through the centuries\/}
\cite[p.~214]{Kl} [emphasis added---authors].
\end{quote}
Klein further points out that such a parallel conception of calculus
\begin{quote}
harks back to old metaphysical speculations concerning the {\em
structure of the continuum\/} according to which this was made up of
[...] infinitely small parts \cite[p.~214]{Kl} [emphasis
added---authors].
\end{quote}

The rival theories of the continuum evoked by Klein can be summarized
as follows.  A Leibnizian definition of the derivative as the
infinitesimal quotient
\[
\frac{\Delta y}{\Delta x},
\] 
whose logical weakness was criticized by Berkeley, was modified by
A.~Robinson by exploiting a map called {\em the standard part\/},
denoted~``st'', from the finite part of a ``thick'' B-continuum (i.e.,
a Bernoullian continuum),%
\footnote{\label{f49}Schubring \cite[p.~170, 173, 187]{Sch} attributes
the first systematic use of infinitesimals as a foundational concept,
to Johann Bernoulli.}
to a ``thin'' A-continuum (i.e., an Archimedean continuum), as
illustrated in Figure~\ref{31}.

This section summarizes a 20th century implementation of an
alternative to an Archimedean continuum, namely an
infinitesimal-enriched continuum.  Such a continuum is not to be
confused with incipient notions of such a continuum found in earlier
centuries.  Johann Bernoulli was one of the first to exploit
infinitesimals in a systematic fashion as a foundational tool in the
calculus.  We will therefore refer to such a continuum as a
Bernoullian continuum, or B-continuum for short.

\begin{figure}
\includegraphics[height=2in]{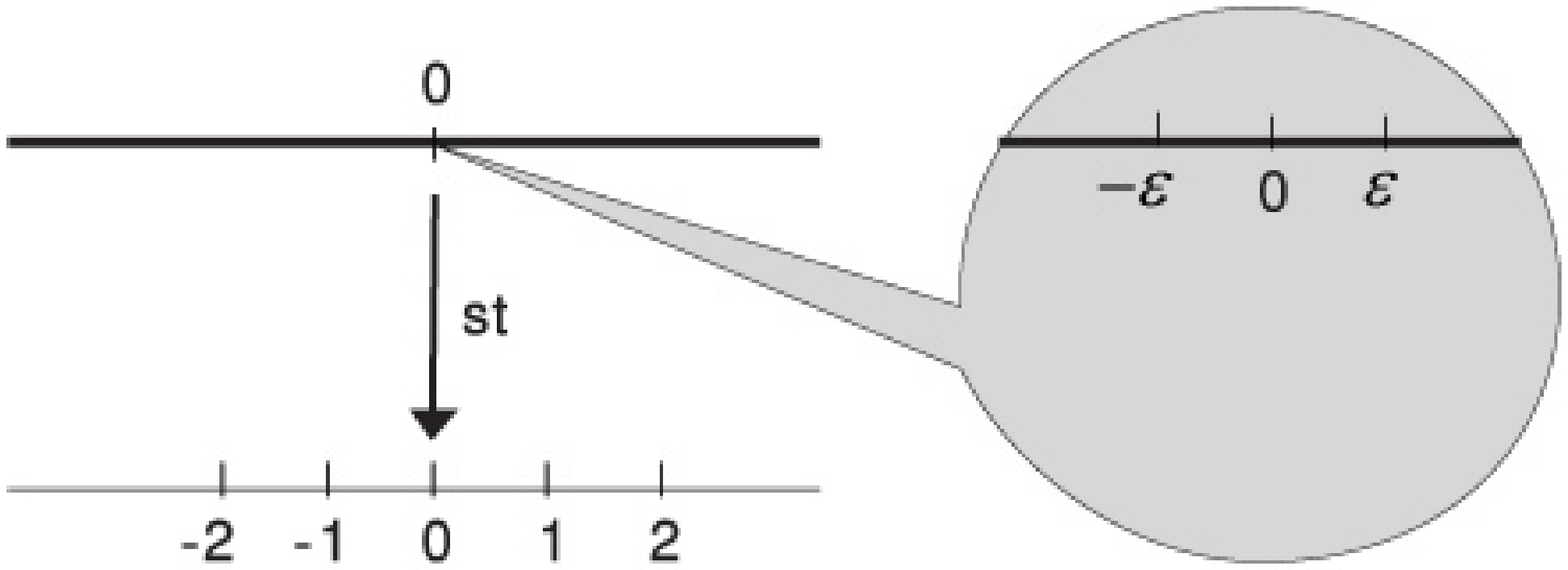}
\caption{\textsf{Zooming in on infinitesimal~$\epsilon$ (here st$(\pm
\epsilon)=0$)}}
\label{tamar}
\end{figure}

\begin{figure}
\includegraphics[height=2.3in]{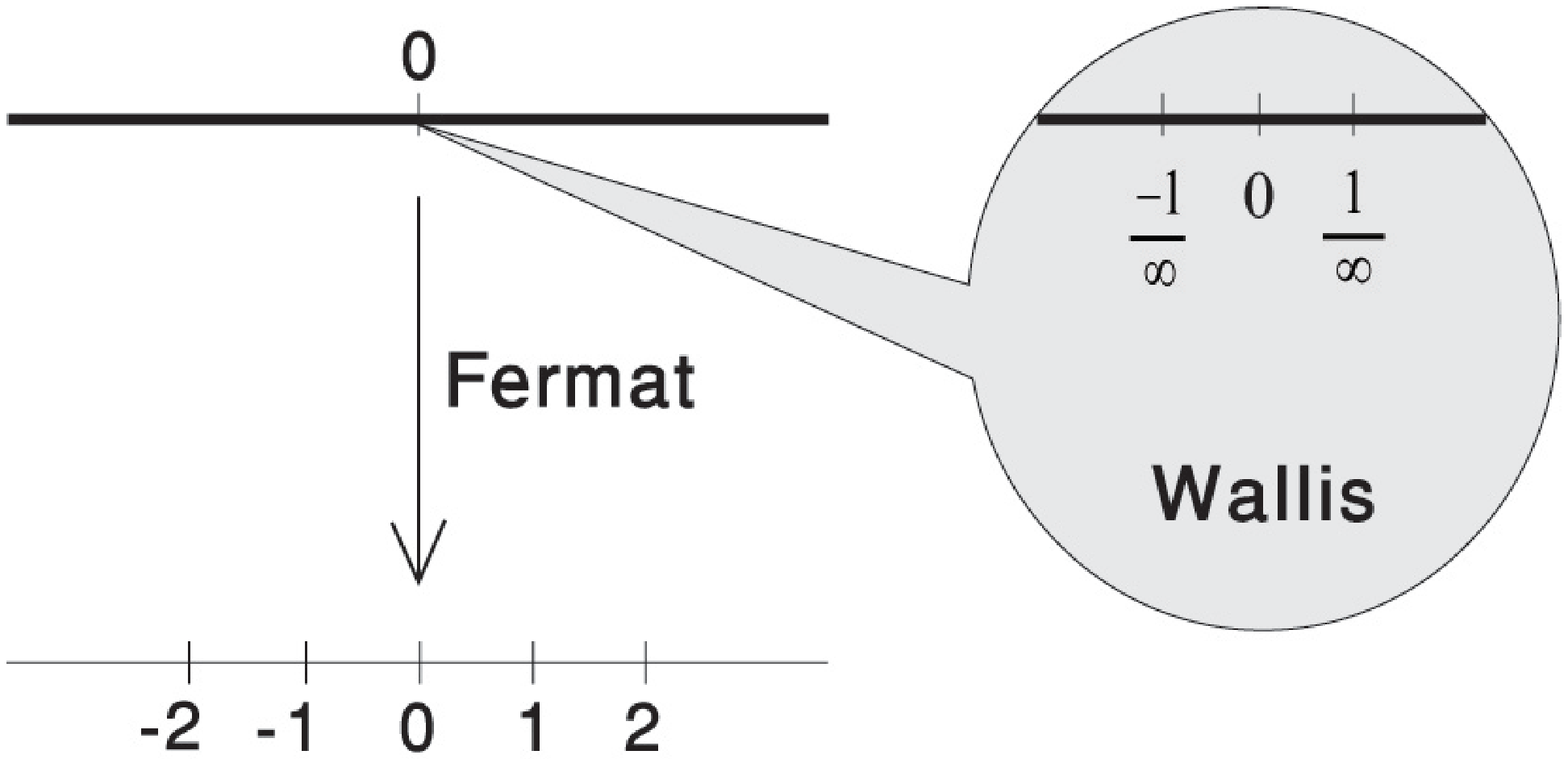}
\caption{\textsf{Zooming in on Wallis's
infinitesimal~$\frac{1}{\infty}$, which is adequal to~$0$ in Fermat's
terminology}}
\label{FermatWallis}
\end{figure}

We illustrate the construction by means of an infinite-resolution
microscope in Figures~\ref{tamar} and~\ref{FermatWallis}.  We will
denote such a B-continuum by the new symbol \RRR{} (``thick-R'').
Such a continuum is constructed in formula~\eqref{bee}.%
\footnote{An alternative implementation of a B-continuum has been
pursued by Lawvere, John L. Bell \cite{Bel08, Bel}, and others.}
We will also denote its finite
part, by
\begin{equation*}
\RRR_{<\infty} = \left\{ x\in \RRR : \; |x|<\infty \right\},
\end{equation*}
so that we have a disjoint union
\begin{equation}
\RRR= \RRR_{<\infty} \cup \RRR_{\infty},
\end{equation}
where~$\RRR_{\infty}$ consists of unlimited hyperreals (i.e., inverses
of nonzero infinitesimals).

The map ``st'' sends each finite point~$x\in \RRR$, to the real point
st$(x)\in \R$ infinitely close to~$x$, see Figure~\ref{tamar}.
Namely, we have:%
\footnote{This is the Fermat-Robinson standard part whose seeds are
found in Fermat's adequality.}
\begin{equation*}
\xymatrix{\quad \RRR_{{<\infty}}^{~} \ar[d]^{{\rm st}} \\ \R}
\end{equation*}
Robinson's answer to Berkeley's {\em logical criticism\/} (see
D.~Sherry \cite{She87}) is to define the derivative as
\begin{equation*}
\hbox{st} \left( \frac{\Delta y}{\Delta x} \right),
\end{equation*}
instead of~$\Delta y/\Delta x$. 

Note that both the term ``hyper-real field'', and an ultrapower
construction thereof, are due to E.~Hewitt in 1948, see
\cite[p.~74]{Hew}.  In 1966, Robinson referred to the 
\begin{quote}
theory of hyperreal fields (Hewitt [1948]) which ... can serve as
non-standard models of analysis \cite[p.~278]{Ro66}.
\end{quote}
The {\em transfer principle\/} is a mathematical implementation of
Leibniz's heuristic {\em law of continuity\/}: ``what succeeds for the
finite numbers succeeds also for the infinite numbers and vice
versa'', see~\cite[p.~266]{Ro66}.  The transfer principle, allowing an
extention of every first-order real statement to the hyperreals, is a
consequence of the theorem of J.~{\L}o{\'s} in 1955, see~\cite{Lo},
and can therefore be referred to as a Leibniz-{\L}o{\'s} transfer
principle.  A Hewitt-{\L}o{\'s} framework allows one to work in a
B-continuum satisfying the transfer principle.  To elaborate on the
ultrapower construction of the hyperreals, let~$\Q^\N$ denote the ring
of sequences of rational numbers.  Let
\begin{equation*}
\left( \Q^\N \right)_C
\end{equation*}
denote the subspace consisting of Cauchy sequences.  The reals are by
definition the quotient field
\begin{equation}
\label{real}
\R:= \left. \left( \Q^\N \right)_C \right/ \mathcal{F}_{\!n\!u\!l\!l},
\end{equation}
where~$\mathcal{F}_{\!n\!u\!l\!l}$ contains all null sequences.
Meanwhile, an infinitesimal-enriched field extension of~$\Q$ may be
obtained by forming the quotient
\begin{equation*}
\left.  \Q^\N \right/ \mathcal{F}_{u}.
\end{equation*}
Here a sequence~$\langle u_n : n\in \N \rangle$ is
in~$\mathcal{F}_{u}$ if and only if the set of indices
\[
\{ n \in \N : u_n = 0 \}
\]
is a member of a fixed ultrafilter.%
\footnote{In this construction, every null sequence defines an
infinitesimal, but the converse is not necessarily true.  Modulo
suitable foundational material, one can ensure that every
infinitesimal is represented by a null sequence; an appropriate
ultrafilter (called a {\em P-point\/}) will exist if one assumes the
continuum hypothesis, or even the weaker Martin's axiom.  See Cutland
{\em et al\/} \cite{CKKR} for details.}
See Figure~\ref{helpful}.

\begin{figure}
\begin{equation*}
\xymatrix{ && \left( \left. \Q^\N \right/ \mathcal{F}_{\!u}
\right)_{<\infty} \ar@{^{(}->} [rr]^{} \ar@{->>}[d]^{\rm st} &&
\RRR_{<\infty} \ar@{->>}[d]^{\rm st} \\ \Q \ar[rr] \ar@{^{(}->} [urr]
&& \R \ar[rr]^{\simeq} && \R }
\end{equation*}
\caption{\textsf{An intermediate field~$\left. \Q^\N \right/
\mathcal{F}_{\!u}$ is built directly out of~$\Q$}}
\label{helpful}
\end{figure}

To give an example, the sequence
\begin{equation}
\label{infinitesimal}
\left\langle \tfrac{(-1)^n}{n} \right\rangle
\end{equation}
represents a nonzero infinitesimal, whose sign depends on whether or
not the set~$2\N$ is a member of the ultrafilter.  To obtain a full
hyperreal field, we replace~$\Q$ by~$\R$ in the construction, and form
a similar quotient
\begin{equation}
\label{bee}
\RRR:= \left.  \R^\N \right/ \mathcal{F}_{u}.
\end{equation}
We wish to emphasize the analogy with formula~\eqref{real} defining
the A-continuum.  Note that, while the leftmost vertical arrow in
Figure~\ref{helpful} is surjective, we have
\begin{equation*}
\left( \Q^\N / \mathcal{F}_{u} \right) \cap \R = \Q.
\end{equation*}
A more detailed discussion of this construction can be found in the
book by M.~Davis~\cite{Dav77}.  
%
%
See also P.~B\l aszczyk \cite{Bl} for some philosophical implications.
More advanced properties of the hyperreals such as saturation were
proved later, see Keisler \cite{Ke94} for a historical outline.  A
helpful ``semicolon'' notation for presenting an extended decimal
expansion of a hyperreal was described by A.~H.~Lightstone~\cite{Li}.
See also P.~Roquette \cite{Roq} for infinitesimal reminiscences.  A
discussion of infinitesimal optics is in K.~Stroyan \cite{Str},
J.~Keisler~\cite{Ke}, D.~Tall~\cite{Ta80}, and L.~Magnani and
R.~Dossena~\cite{MD, DM}.  

Applications of the B-continuum range from aid in teaching calculus
\cite{El, KK1, KK2, Ta91, Ta09a} to the Bolzmann equation (see
L.~Arkeryd~\cite{Ar81, Ar05}); modeling of timed systems in computer
science (see H.~Rust \cite{Ru}); mathematical economics (see Anderson
\cite{An00}); mathematical physics (see Albeverio {\em et al.\/}
\cite{Alb}); etc.

\section*{Acknowledgments}

The authors are grateful to Piotr B\l aszczyk, Paolo Giordano, Thomas
Mormann, and David Tall for valuable comments that helped improve the
text.  Hilton Kramer's influence is obvious throughout.

\end{document}